\documentclass[12pt]{article}%
\usepackage{amsfonts}
\usepackage{sw20bams}
\usepackage{amsmath}
\usepackage{amssymb}
\usepackage{graphicx}%
\providecommand{\U}[1]{\protect\rule{.1in}{.1in}}
\begin{document}

\title{Geo/Geo/$2$ Queues and the Poisson Clumping Heuristic}
\author{Steven Finch}
\date{February 25, 2019}
\maketitle

\begin{abstract}
In discrete time, customers arrive at random. \ Each waits until one of two
servers is available; each thereafter departs at random. \ We seek the
distribution of maximum line length of idle customers. \ In the context of an
emergency room (for medical treatment), the virtue of one fast doctor over two
slow doctors is explored. \ Via limiting argument to continuous time, we study
likewise the M/M/$2$ queue.

\end{abstract}

\footnotetext{Copyright \copyright \ 2019 by Steven R. Finch. All rights
reserved.}Let $0<r<1$ and $0<p<r$. \ Consider the Julia program:%
\[%
\begin{array}
[c]{l}%
\text{\texttt{u = 0}}\\
\text{\texttt{m = 0}}\\
\text{\texttt{for t=1:n}}\\
\text{\texttt{\ \ x = rand()%
$<$%
p \ \ \ \ \ \ \ \ \ \ \# x=1 means that an arrival occurs}}\\
\text{\texttt{\ \ y = rand()%
$<$%
r \ \ \ \ \ \ \ \ \ \ \# y=1 means that one departure occurs}}\\
\text{\texttt{\ \ if u==0 }}\\
\text{\texttt{\ \ \ \ u = x\ \ \ \ \ \ \ \ \ \ \ \ \ \ \ \ \ \# increment is 1
or 0\ \ }}\\
\text{\texttt{\ \ else}}\\
\text{\texttt{\ \ \ \ u = max(0,u+x-y) \ \ \ \ \ \# increment is 1, 0 or
-1\ }}\\
\text{\texttt{\ \ end}}\\
\text{\texttt{\ \ m = max(m,u)}}\\
\text{\texttt{end}}\\
\text{\texttt{return m}}%
\end{array}
\]
which simulates the maximum value of a Geo/Geo/$1$ queue with LAS-DA over $n$
time steps. \ The Boolean expressions containing Julia's Uniform $[0,1]$
random deviate generator ensure that $X\sim\,$Bernoulli($p$) and $Y\sim
\,$Bernoulli($r$). \ The word \textquotedblleft Geometric\textquotedblright%
\ arises because%
\[%
\begin{array}
[c]{ccc}%
\mathbb{P}\left\{  \text{time lapse between adjacent arrivals is }i\right\}
=p\,q^{i-1}, &  & i\geq1
\end{array}
\]
where $q=1-p$ and
\[%
\begin{array}
[c]{ccc}%
\mathbb{P}\left\{  \text{time lapse between adjacent departures is }j\right\}
=r\,s^{j-1}, &  & j\geq1
\end{array}
\]
where $s=1-r$. \ Clearly $s<q<1$. \ LAS\ stands for \textquotedblleft late
arrival system\textquotedblright\ and DA\ stands for \textquotedblleft delayed
access\textquotedblright\ \cite{Htr-heu}; in particular, a customer entering
an empty queue at time $t$ is not immediately eligible for service, but rather
at time $t+1$. We study the asymptotic distribution of the maximum $M_{n}$ in
Section 1.

Now let $0<r<1$ and $0<p<2r$. \ Consider the program:
\[%
\begin{array}
[c]{l}%
\text{\texttt{u = 0}}\\
\text{\texttt{m = 0}}\\
\text{\texttt{for t=1:n}}\\
\text{\texttt{\ \ x = rand()%
$<$%
p \ \ \ \ \ \ \ \ \ \ \# x=1 means that an arrival occurs}}\\
\text{\texttt{\ \ y1 = rand()%
$<$%
r \ \ \ \ \ \ \ \ \ \# y1+y2=2 means that two departures occur}}\\
\text{\texttt{\ \ y2 = rand()%
$<$%
r \ \ \ \ \ }}\\
\text{\texttt{\ \ if u==0 }}\\
\text{\texttt{\ \ \ \ u = x\ \ \ \ \ \ \ \ \ \ \ \ \ \ \ \ \ \ \ \ \ \ \ \ \#
increment is 1 or 0\ \ }}\\
\text{\texttt{\ \ else}}\\
\text{\texttt{\ \ \ \ if u==1}}\\
\text{\texttt{\ \ \ \ \ \ u = max(0,u+x-y1) \ \ \ \ \ \ \ \ \ \# increment is
1, 0 or -1}}\\
\text{\texttt{\ \ \ \ else}}\\
\text{\texttt{\ \ \ \ \ \ u = max(0,u+x-y1-y2)\ \ \ \ \ \ \ \# increment is 1,
0, -1 or -2}}\\
\text{\texttt{\ \ \ \ end}}\\
\text{\texttt{\ \ end}}\\
\text{\texttt{\ \ m = max(m,u)}}\\
\text{\texttt{end}}\\
\text{\texttt{return m}}%
\end{array}
\]
which simulates the maximum value of a Geo/Geo/$2$ queue with LAS-DA over $n$
time steps. \ Clearly again $2s-1<q<1$. We study the asymptotic distribution
of the maximum $M_{n}$ in Section 2. \ Taking time steps to be not of length
$1$, but of length $\Delta>0$ and allowing $\Delta$ to approach $0$, gives
analogous formulas for the M/M/$2$ queue in Section 3.

The Poisson clumping heuristic \cite{Ald-heu}, while not a theorem, gives
results identical to exact asymptotic expressions when such exist, and
evidently provides excellent predictions otherwise. \ Consider an irreducible
positive recurrent Markov chain with stationary distribution $\pi$. \ For
large enough $k$, the maximum of the chain satisfies%
\[
\mathbb{P}\left\{  M_{n}<k\right\}  \sim\exp\left(  -\frac{\pi_{k}}%
{\mathbb{E}(C)}n\right)
\]
as $n\rightarrow\infty$, where $C$ is the sojourn time in $k$ during a clump
of nearby visits to $k$. \ 

\section{One Server}

Starting with transition matrix%
\[
\left(
\begin{array}
[c]{cccccc}%
q & p & 0 & 0 & 0 & \cdots\\
qr & pr+qs & ps & 0 & 0 & \cdots\\
0 & qr & pr+qs & ps & 0 & \cdots\\
0 & 0 & qr & pr+qs & ps & \cdots\\
0 & 0 & 0 & qr & pr+qs & \cdots\\
\vdots & \vdots & \vdots & \vdots & \vdots & \ddots
\end{array}
\right)
\]
we obtain \cite{AH-heu, Alf-heu}%
\[
\pi_{0}=\frac{r-p}{r}=\frac{qr}{p}\pi_{1},
\]%
\[%
\begin{array}
[c]{ccccc}%
\pi_{j}=\omega^{j-1}\pi_{1} &  & \text{for} &  & j\geq2
\end{array}
\]
where $0<\omega<1$ satisfies the quadratic equation $\omega=(q\omega
+p)(r\omega+s)$, that is,%
\[
\omega=\frac{ps}{qr}
\]
and%
\[
\pi_{1}=\frac{p(1-\omega)}{r}.
\]
Note that, if $k=\log_{1/\omega}(n)+h+1$, we have%
\[
\left(  \dfrac{1}{\omega}\right)  ^{k}=n\left(  \dfrac{1}{\omega}\right)
^{h+1}
\]
thus%
\[
\pi_{k}n=\frac{p(1-\omega)}{r}\omega^{k-1}n=\frac{p(1-\omega)}{r}\omega
^{h}=\frac{p(r-p)}{qr^{2}}\left(  \frac{ps}{qr}\right)  ^{h}.
\]
Following the argument in section 1.3 of \cite{Fi0-heu}, replacing $p^{2}$ by
$ps$, $q^{2}$ by $qr$ and $2pq$ by $pr+qs$, we obtain%
\[
\mathbb{E}(C)=\frac{1}{1-(p+s)}=\frac{1}{r-p}
\]
which implies%
\begin{align*}
\mathbb{P}\left\{  M_{n}\leq\log_{1/\omega}(n)+h\right\}   &  =P\left\{
M_{n}<\log_{1/\omega}(n)+h+1\right\} \\
&  \sim\exp\left[  -\frac{p(r-p)^{2}}{qr^{2}}\left(  \frac{ps}{qr}\right)
^{h}\right]
\end{align*}
as $n\rightarrow\infty$. \ For example, if%
\[%
\begin{array}
[c]{ccccccc}%
p=\frac{1}{3}, &  & q=\frac{2}{3}, &  & r=\frac{1}{2}, &  & s=\frac{1}{2}%
\end{array}
\]
we have
\begin{align*}
\mathbb{E}\left(  M_{n}\right)   &  \approx\frac{\ln(n)}{\ln(\frac{qr}{ps}%
)}+\frac{\gamma+\ln\left(  \frac{p(r-p)^{2}}{qr^{2}}\right)  }{\ln(\frac
{qr}{ps})}+\frac{1}{2}\\
&  \approx(1.4426950408...)\ln(n)-(2.8371788241...)
\end{align*}
for sufficiently large $n$, where $\gamma$ denotes Euler's constant
\cite{Fi1-heu}.

\section{Two Servers}

Starting with transition matrix%
\[
\left(
\begin{array}
[c]{ccccccc}%
q & p & 0 & 0 & 0 & 0 & \cdots\\
qr & pr+qs & ps & 0 & 0 & 0 & \cdots\\
qr^{2} & pr^{2}+2qrs & 2prs+qs^{2} & ps^{2} & 0 & 0 & \cdots\\
0 & qr^{2} & pr^{2}+2qrs & 2prs+qs^{2} & ps^{2} & 0 & \cdots\\
0 & 0 & qr^{2} & pr^{2}+2qrs & 2prs+qs^{2} & ps^{2} & \cdots\\
0 & 0 & 0 & qr^{2} & pr^{2}+2qrs & 2prs+qs^{2} & \cdots\\
\vdots & \vdots & \vdots & \vdots & \vdots & \vdots & \ddots
\end{array}
\right)
\]
we obtain \cite{AH-heu, NPR-heu}%
\[
\pi_{0}=\frac{qr^{2}}{p^{2}s}\left(  1+qs+qr\omega\right)  \pi_{2},
\]%
\[
\pi_{1}=\frac{r}{ps}\left(  r+2qs+qr\omega\right)  \pi_{2},
\]%
\[%
\begin{array}
[c]{ccccc}%
\pi_{j}=\omega^{j-2}\pi_{2} &  & \text{for} &  & j\geq3
\end{array}
\]
where $0<\omega<1$ satisfies the cubic equation $\omega=(q\omega
+p)(r\omega+s)^{2}$, that is,%
\[%
\begin{array}
[c]{ccc}%
\omega=\dfrac{-r-2qs+\theta}{2qr}, &  & \theta=\sqrt{r^{2}+4qs}%
\end{array}
\]
and%
\[
\pi_{2}=\frac{p^{2}s(1-\omega)}{p^{2}s+r\left[  r+pqs+q(p+qr)(s+r\omega
)\right]  (1-\omega)}.
\]
Note that, if $k=\log_{1/\omega}(n)+h+1$, we have%
\[
\left(  \dfrac{1}{\omega}\right)  ^{k}=n\left(  \dfrac{1}{\omega}\right)
^{h+1}
\]
thus%
\[
\pi_{k}n=\pi_{2}\omega^{k-2}n=\pi_{2}\omega^{h-1}.
\]

We need now to calculate $\mathbb{E}(C)$. Consider a random walk on the
integers consisting of incremental steps satisfying%
\[
\left\{
\begin{array}
[c]{lll}%
-2 &  & \text{with probability }qr^{2},\\
-1 &  & \text{with probability }pr^{2}+2qrs,\\
0 &  & \text{with probability }2prs+qs^{2},\\
1 &  & \text{with probability }ps^{2}.
\end{array}
\right.
\]
For nonzero $j$, let $\nu_{j}$ denote the probability that, starting from
$-j$, the walker eventually hits $0$. \ Let $\nu_{0}$ denote the probability
that, starting from $0$, the walker eventually returns to $0$ (at some future
time). \ We have two values for $\nu_{0}$: when it is used in a recursion, it
is equal to $1$; when it corresponds to a return probability, it retains the
symbol $\nu_{0}$. \ Using%
\[%
\begin{array}
[c]{ccc}%
\nu_{j}=ps^{2}\nu_{j-1}+(2prs+qs^{2})\nu_{j}+(pr^{2}+2qrs)\nu_{j+1}+qr^{2}%
\nu_{j+2}, &  & j\geq1;
\end{array}
\]%
\[
\nu_{0}=ps^{2}\nu_{-1}+(2prs+qs^{2})+(pr^{2}+2qrs)\nu_{1}+qr^{2}\nu_{2}
\]
define%
\begin{align*}
F(z)  &  =%
{\displaystyle\sum\limits_{j=1}^{\infty}}
\nu_{j}z^{j}\\
&  =ps^{2}z%
{\displaystyle\sum\limits_{j=1}^{\infty}}
\nu_{j-1}z^{j-1}+(2prs+qs^{2})%
{\displaystyle\sum\limits_{j=1}^{\infty}}
\nu_{j}z^{j}+\frac{pr^{2}+2qrs}{z}%
{\displaystyle\sum\limits_{j=1}^{\infty}}
\nu_{j+1}z^{j+1}\\
&  +\frac{qr^{2}}{z^{2}}%
{\displaystyle\sum\limits_{j=1}^{\infty}}
\nu_{j+2}z^{j+2}\\
&  =ps^{2}z\left[  F(z)+1\right]  +(2prs+qs^{2})F(z)+\frac{pr^{2}+2qrs}%
{z}\left[  F(z)-\nu_{1}z\right] \\
&  +\frac{qr^{2}}{z^{2}}\left[  F(z)-\nu_{1}z-\nu_{2}z^{2}\right]
\end{align*}
equivalently%
\begin{align*}
&  \left[  1-ps^{2}z-2prs-qs^{2}-\frac{pr^{2}+2qrs}{z}-\frac{qr^{2}}{z^{2}%
}\right]  F(z)\\
&  =ps^{2}z-\frac{pr^{2}+2qrs}{z}(\nu_{1}z)-\frac{qr^{2}}{z^{2}}(\nu_{1}%
z+\nu_{2}z^{2})
\end{align*}
equivalently%
\begin{align*}
&  \left[  qr^{2}+(pr^{2}+2qrs)z-(1-2prs-qs^{2})z^{2}+ps^{2}z^{3}\right]
F(z)\\
&  =-ps^{2}z^{3}+(pr^{2}+2qrs)z(\nu_{1}z)+qr^{2}(\nu_{1}z+\nu_{2}z^{2})
\end{align*}
equivalently%
\begin{align*}
&  (1-z)\left[  qr^{2}+(2qs+r)rz-ps^{2}z^{2}\right]  F(z)\\
&  =-ps^{2}z^{3}+pr^{2}z^{2}\nu_{1}+2qrsz^{2}\nu_{1}+qr^{2}z\nu_{1}\\
&  +z^{2}\left(  \nu_{0}-ps^{2}\nu_{-1}-2prs-qs^{2}-pr^{2}\nu_{1}-2qrs\nu
_{1}\right) \\
&  =z^{2}\nu_{0}+qr^{2}z\nu_{1}-ps^{2}z^{2}\nu_{-1}-2prsz^{2}-qs^{2}%
z^{2}-ps^{2}z^{3}.
\end{align*}
Examine the denominator of $F(z)$. \ Only the first two of its three zeroes
$z_{1}$, $1$, $z_{2}$ are of interest (the third is $>1$). \ Note that%
\[
z_{1}=\frac{(2qs+r-\theta)r}{2ps^{2}}.
\]
Substituting $z=1$ and $z=z_{1}$ into the numerator $N_{F}$ of $F(z)$, then
setting $N_{F}=0$, gives two equations in three unknowns. \ At this point in
section 1.3 of \cite{Fi0-heu}, we utilized a simple formula for $\nu_{-1}$ in
terms of $\nu_{1}$. \ Due to the complexity of $F(z)$ here, a different
approach is required. \ 

Using%
\[%
\begin{array}
[c]{ccc}%
\nu_{-j}=ps^{2}\nu_{-j-1}+(2prs+qs^{2})\nu_{-j}+(pr^{2}+2qrs)\nu_{-j+1}%
+qr^{2}\nu_{-j+2}, &  & j\geq1;
\end{array}
\]%
\[
\nu_{0}=ps^{2}\nu_{-1}+(2prs+qs^{2})+(pr^{2}+2qrs)\nu_{1}+qr^{2}\nu_{2}
\]
define%
\begin{align*}
G(z)  &  =%
{\displaystyle\sum\limits_{j=1}^{\infty}}
\nu_{-j}z^{j}\\
&  =\frac{ps^{2}}{z}%
{\displaystyle\sum\limits_{j=1}^{\infty}}
\nu_{-j-1}z^{j+1}+(2prs+qs^{2})%
{\displaystyle\sum\limits_{j=1}^{\infty}}
\nu_{-j}z^{j}+(pr^{2}+2qrs)z%
{\displaystyle\sum\limits_{j=1}^{\infty}}
\nu_{-j+1}z^{j-1}\\
&  +qr^{2}z^{2}%
{\displaystyle\sum\limits_{j=1}^{\infty}}
\nu_{-j+2}z^{j-2}\\
&  =\frac{ps^{2}}{z}\left[  G(z)-\nu_{-1}z\right]  +(2prs+qs^{2}%
)G(z)+(pr^{2}+2qrs)z\left[  G(z)+1\right] \\
&  +qr^{2}z^{2}\left[  G(z)+1+\frac{\nu_{1}}{z}\right]
\end{align*}
equivalently%
\begin{align*}
&  \left[  1-\frac{ps^{2}}{z}-2prs-qs^{2}-(pr^{2}+2qrs)z-qr^{2}z^{2}\right]
G(z)\\
&  =-\frac{ps^{2}}{z}(\nu_{-1}z)+(pr^{2}+2qrs)z+qr^{2}z^{2}\left(  1+\frac
{\nu_{1}}{z}\right)
\end{align*}
equivalently%
\begin{align*}
&  \left[  qr^{2}z^{3}+(pr^{2}+2qrs)z^{2}-(1-2prs-qs^{2})z+ps^{2}\right]
G(z)\\
&  =ps^{2}(\nu_{-1}z)-(pr^{2}+2qrs)z^{2}-qr^{2}z^{2}\left(  z+\nu_{1}\right)
\end{align*}
equivalently%
\begin{align*}
&  (1-z)\left[  ps^{2}-(2qs+r)rz-qr^{2}z^{2}\right]  G(z)\\
&  =ps^{2}z\nu_{-1}-pr^{2}z^{2}-2qrsz^{2}-qr^{2}z^{3}-qr^{2}z^{2}\nu_{1}.
\end{align*}
Examine the denominator of $G(z)$. \ Only the zero $z_{3}$ of smallest
modulus:%
\[
z_{3}=\frac{-(2qs+r-\theta)r}{2qr^{2}}=\frac{-2qs-r+\theta}{2qr}
\]
interests us. \ Substituting $z=z_{3}$ into the numerator $N_{G}$ of $G(z)$,
and setting $N_{G}=0$, gives a third equation (to include with the other two
from earlier). \ Solving the simultaneous system in $\nu_{0}$, $\nu_{-1}$,
$\nu_{1}$, we obtain%

\[
\nu_{0}=\frac{6q-4qr+r^{2}-2q\theta-r\theta}{2q},
\]%
\[
\nu_{-1}=\frac{(2qs+r-\theta)(qr-\theta)}{2pqs^{2}},
\]%
\[
\nu_{1}=\frac{-r-2qs+\theta}{2qr}=\omega
\]
which implies%
\begin{align*}
\mathbb{P}\left\{  M_{n}\leq\log_{1/\omega}(n)+h\right\}   &  =P\left\{
M_{n}<\log_{1/\omega}(n)+h+1\right\} \\
&  \sim\exp\left[  -\frac{\pi_{2}(1-\nu_{0})}{\omega}\omega^{h}\right]
\end{align*}
as $n\rightarrow\infty$. \ For example, if%
\[%
\begin{array}
[c]{ccccccc}%
p=\frac{1}{3}, &  & q=\frac{2}{3}, &  & r=\frac{1}{4}, &  & s=\frac{3}{4}%
\end{array}
\]
we have%
\[%
\begin{array}
[c]{ccc}%
\omega=0.5584219849..., &  & \pi_{2}=0.2270554252...,
\end{array}
\]%
\[%
\begin{array}
[c]{ccc}%
\nu_{0}=0.8414579643..., &  & \dfrac{\pi_{2}(1-\nu_{0})}{\omega}%
=0.0644634887...,
\end{array}
\]%
\begin{align*}
\mathbb{E}\left(  M_{n}\right)   &  \approx\frac{\ln(n)}{\ln(\frac{1}{\omega
})}+\frac{\gamma+\ln\left(  \frac{\pi_{2}(1-\nu_{0})}{\omega}\right)  }%
{\ln(\frac{1}{\omega})}+\frac{1}{2}\\
&  \approx(1.7163246381...)\ln(n)-(3.2148827577...)
\end{align*}
for sufficiently large $n$.

The use of an expected maximum for performance analysis, instead of a simple
average, does not appear to lead to surprising outcomes. \ A\ corollary of the
preceding numerical results is that, in a busy hospital emergency room (with
$p=1/3$), one fast doctor (with $r=1/2$) outperforms two slow doctors (each
with $r=1/4$). \ For average queue lengths \cite{AH-heu},%
\[%
{\displaystyle\sum\limits_{j=1}^{\infty}}
j\pi_{j}=\frac{1}{(1-\omega)^{2}}\pi_{1}=\frac{pq}{r-p}=1.33333...
\]
corresponding to Geo/Geo/$1$ and
\[%
{\displaystyle\sum\limits_{j=1}^{\infty}}
j\pi_{j}=\pi_{1}+\frac{2-\omega}{(1-\omega)^{2}}\pi_{2}=1.98358...
\]
corresponding to Geo/Geo/$2$. \ This is also consistent with results in
\cite{Fi2-heu} governing deterministic traffic signals: we do better with an
$RGRG...$ pattern than with $RRGG...$.

\section{From Discrete to Continuous}

Consider an M/M/$1$ queue with arrival rate $\lambda$ and service rate $\mu$.
\ If $\lambda<\mu$, then parameters of a Geo/Geo/$1$ queue with $p=\lambda
\Delta$ and $r=\mu\Delta$ approach those of the M/M/$1$ queue as
$\Delta\rightarrow0^{+}$. In particular \cite{AH-heu},%
\[
\lim_{\Delta\rightarrow0^{+}}\pi_{k}n=\lim_{\Delta\rightarrow0^{+}}%
\frac{p(r-p)}{qr^{2}}\left(  \frac{ps}{qr}\right)  ^{h}=\frac{\mu-\lambda}%
{\mu}\left(  \frac{\lambda}{\mu}\right)  ^{h+1},
\]

\[
\lim_{\Delta\rightarrow0^{+}}\frac{1}{\mathbb{E}(C)\Delta}=\lim_{\Delta
\rightarrow0^{+}}\frac{r-p}{\Delta}=\mu-\lambda
\]
and hence, over the time interval $[0,x]$,
\[
\mathbb{P}\left\{  M_{x}\leq\log_{\mu/\lambda}(x)+h\right\}  \sim\exp\left[
-\frac{(\mu-\lambda)^{2}}{\mu}\left(  \frac{\lambda}{\mu}\right)
^{h+2}\right]
\]
as $x\rightarrow\infty$, consistent with \cite{Ald-heu}. \ For $\lambda=1/3$
and $\mu=1/2$, we have%
\begin{align*}
\mathbb{E}\left(  M_{x}\right)   &  \approx\frac{\ln(x)}{\ln(\frac{\mu
}{\lambda})}+\frac{\gamma+\ln\left(  \frac{\lambda^{2}(\mu-\lambda)^{2}}%
{\mu^{3}}\right)  }{\ln(\frac{\mu}{\lambda})}+\frac{1}{2}\\
&  \approx(2.4663034623...)\ln(x)-(7.2049448811...).
\end{align*}
Consider instead an M/M/$2$ queue with arrival rate $\lambda$ and service rate
$\mu$. \ If $\lambda<2\mu$, then \cite{AH-heu}%
\[
\lim_{\Delta\rightarrow0^{+}}\pi_{k}n=\lim_{\Delta\rightarrow0^{+}}\frac
{\pi_{2}}{\omega}\omega^{h}=2\frac{2\mu-\lambda}{2\mu+\lambda}\left(
\frac{\lambda}{2\mu}\right)  ^{h+1},
\]

\[
\lim_{\Delta\rightarrow0^{+}}\frac{1-\nu_{0}}{\Delta}=\lim_{\Delta
\rightarrow0^{+}}\frac{1}{\Delta}\left(  \frac{6q-4qr+r^{2}-2q\theta-r\theta
}{2q}\right)  =2\mu-\lambda
\]
and hence, over the time interval $[0,x]$,%
\[
\mathbb{P}\left\{  M_{x}\leq\log_{2\mu/\lambda}(x)+h\right\}  \sim\exp\left[
-2\frac{(2\mu-\lambda)^{2}}{2\mu+\lambda}\left(  \frac{\lambda}{2\mu}\right)
^{h+2}\right]
\]
as $x\rightarrow\infty$. \ A reference for this formula is not known. \ For
$\lambda=1/3$ and $\mu=1/4$, we have
\begin{align*}
\mathbb{E}\left(  M_{x}\right)   &  \approx\frac{\ln(x)}{\ln(\frac{2\mu
}{\lambda})}+\frac{\gamma+\ln\left(  \frac{\lambda^{2}(2\mu-\lambda)^{2}}%
{2\mu^{2}(2\mu+\lambda)}\right)  }{\ln(\frac{2\mu}{\lambda})}+\frac{1}{2}\\
&  \approx(2.4663034623...)\ln(x)-(6.7552845943...).
\end{align*}
Again, with regard to expected maximums, in an emergency room ($\lambda=1/3$),
one fast doctor ($\mu=1/2$) outperforms two slow doctors (each $\mu=1/4$).
\ Well-known formulas for simple averages \cite{AH-heu, HL-heu} are instead%
\[
\frac{1}{(1-\frac{\lambda}{\mu})^{2}}\left(  \lim_{\Delta\rightarrow0^{+}}%
\pi_{1}\right)  =\frac{\lambda}{\mu-\lambda}=2
\]
corresponding to M/M/$1$ and
\[
\left(  \lim_{\Delta\rightarrow0^{+}}\pi_{1}\right)  +\frac{2-\frac{\lambda
}{2\mu}}{(1-\frac{\lambda}{2\mu})^{2}}\left(  \lim_{\Delta\rightarrow0^{+}}%
\pi_{2}\right)  =\frac{4\lambda\mu}{\left(  2\mu-\lambda\right)  \left(
2\mu+\lambda\right)  }=2.4
\]
corresponding to M/M/$2$. \ Results summarizing a\ continuous-time analog of
deterministic traffic signals would be good to see someday.

\section{Appendix I}

Let $0<\lambda<c\mu$. \ Consider the R program:%
\[%
\begin{array}
[c]{l}%
\text{\texttt{K
$<$%
- rpois(1,x*lambda)}}\\
\text{\texttt{P
$<$%
- matrix(0,K,3) \ \ \ \ \ \ \ \ \# matrix of patients}}\\
\text{\texttt{P[,1]
$<$%
- sort(runif(K,0,x)) \ \ \ \ \ \ \ \ \# arrival times}}\\
\text{\texttt{P[,3]
$<$%
- rexp(K,mu) \ \ \ \ \ \ \ \ \ \ \ \ \ \ \ \ \# treatment lengths}}\\
\text{\texttt{D
$<$%
- rep(0,c) \ \ \ \ \ \ \ \ \ \ \ \ \ \# vector of doctors}}\\
\text{\texttt{k.sys
$<$%
- function(i,P) length(P[P[,1]%
$<$%
P[i,1] \& P[i,1]%
$<$%
P[,2]+P[,3],1])}}\\
\text{\texttt{k.que
$<$%
- function(i,P) length(P[P[,1]%
$<$%
P[i,1] \& P[i,1]%
$<$%
P[,2],1])}}\\
\text{\texttt{for (i in 1:K)}}\\
\text{\texttt{\ \ \{}}\\
\text{\texttt{\ \ \ \ j
$<$%
- which.min(D)}}\\
\text{\texttt{\ \ \ \ P[i,2]
$<$%
- max(P[i,1],D[j])}}\\
\text{\texttt{\ \ \ \ D[j]
$<$%
- P[i,2] + P[i,3] \ \ \ \ \ \ \ \ \# departure times}}\\
\text{\texttt{\ \ \}}}\\
\text{\texttt{L.sys
$<$%
- sapply(1:K,k.sys,P=P)}}\\
\text{\texttt{L.que
$<$%
- sapply(1:K,k.que,P=P)}}\\
\text{\texttt{list(max(L.sys),max(L.que))}}%
\end{array}
\]
which simulates the maximum value of an M/M/$c$ queue over the time interval
$[0,x]$. \ More precisely, at any arrival time $t=P_{i,1}$, let $L_{sys}$
denote the number of patients in the \textbf{system} (either queue or
treatment) and $L_{que}$ denote the number of patients in the \textbf{queue}
(excluding treatment). \ The maximums of $L_{sys}$ and $L_{que}$ over all
arrival times up to $x$ satisfy%
\[%
\begin{array}
[c]{ccc}%
\max\limits_{0\leq t\leq x}L_{sys}=c+\max\limits_{0\leq t\leq x}L_{que} &  &
\text{almost always,}%
\end{array}
\]
for large enough $n$. \ Simulation further suggests that $\max\nolimits_{0\leq
t\leq x}L_{sys}$ possesses the same distribution as $M_{x}$ defined in Section
3. \ This is somewhat surprising because $M_{x}$ is the limit (as
$\Delta\rightarrow0^{+}$) of $M_{n}$ which, in turn, is based \emph{not} on
Geo/Geo/$c$ system lengths but rather queue lengths. \ A\ resolution of this
minor mystery would be welcome.

\section{Appendix II}

For clarity's sake, consider the (simplified) Julia program:%
\[%
\begin{array}
[c]{l}%
\text{\texttt{u = 0}}\\
\text{\texttt{m = 0}}\\
\text{\texttt{for t=1:n}}\\
\text{\texttt{\ \ x = rand()%
$<$%
p \ \ \ \ \ \ \ \ \ \ \# x=1 means that an arrival occurs}}\\
\text{\texttt{\ \ y = rand()%
$<$%
r \ \ \ \ \ \ \ \ \ \ \# y=1 means that one departure occurs}}\\
\text{\texttt{\ \ u = max(0,u+x-y) \ \ \ \ \ \ \ \# increment is 1, 0 or
-1\ }}\\
\text{\texttt{\ \ m = max(m,u)}}\\
\text{\texttt{end}}\\
\text{\texttt{return m}}%
\end{array}
\]
where $0<p<r<1$. \ The transition matrix for this, a Geo/Geo/$1$ queue with
EAS (\textquotedblleft early arrival system\textquotedblright), is%
\[
\left(
\begin{array}
[c]{cccccc}%
pr+q & ps & 0 & 0 & 0 & \cdots\\
qr & pr+qs & ps & 0 & 0 & \cdots\\
0 & qr & pr+qs & ps & 0 & \cdots\\
0 & 0 & qr & pr+qs & ps & \cdots\\
0 & 0 & 0 & qr & pr+qs & \cdots\\
\vdots & \vdots & \vdots & \vdots & \vdots & \ddots
\end{array}
\right)  .
\]
Guided by reasoning in section 1.1 of \cite{Fi0-heu}, substituting $p^{2}$ by
$ps$, $q^{2}$ by $qr$ and $2pq$ by $pr+qs$, we obtain
\[%
\begin{array}
[c]{ccccc}%
\pi_{j}=(1-\omega)\omega^{j} &  & \text{for} &  & j\geq0
\end{array}
\]
where $\omega=(ps)/(qr)$ as before. $\ $From $k=\log_{1/\omega}(n)+h+1$
follows%
\[
\left(  \dfrac{1}{\omega}\right)  ^{k}=n\left(  \dfrac{1}{\omega}\right)
^{h+1}%
\]
thus%
\[
\pi_{k}n=(1-\omega)\omega^{k}n=(1-\omega)\omega^{h+1}=\frac{r-p}{qr}\left(
\frac{ps}{qr}\right)  ^{h+1}.
\]
The clumping heuristic, coupled with $\mathbb{E}(C)=1/(r-p)$, guarantees%
\[
\mathbb{P}\left\{  M_{n}\leq\log_{1/\omega}(n)+h\right\}  \sim\exp\left[
-\frac{ps(r-p)^{2}}{q^{2}r^{2}}\left(  \frac{ps}{qr}\right)  ^{h}\right]  ,
\]%
\[
\mathbb{E}\left(  M_{n}\right)  \approx\frac{\ln(n)}{\ln(\frac{qr}{ps})}%
+\frac{\gamma+\ln\left(  \frac{ps(r-p)^{2}}{q^{2}r^{2}}\right)  }{\ln
(\frac{qr}{ps})}+\frac{1}{2}%
\]
as $n\rightarrow\infty$. \ It is natural to question whether there exists a
derivation of such formulas that does not depend on the truth of an unproven assertion.

One possible answer is to imagine the EAS increments $Z_{1}$, $Z_{2}$, \ldots,
$Z_{n}$ as a lazy random walk:%
\[%
\begin{array}
[c]{ccccc}%
\mathbb{P\{}Z_{t}=1\}=a, &  & \mathbb{P\{}Z_{t}=-1\}=b, &  & \mathbb{P\{}%
Z_{t}=0\}=c
\end{array}
\]
with reflection at the origin, giving \cite{Fi3-heu}%
\[
\mathbb{E}^{\prime}\left(  M_{n}\right)  \approx\frac{\ln\left(
(a+b)n\right)  }{\ln(\frac{b}{a})}+\frac{\gamma+\ln\left(  \frac{a(b-a)^{2}%
}{b^{2}}\right)  }{\ln(\frac{b}{a})}+\frac{1}{2}.
\]
(Reason: the laziness effectively reduces the sample size by a factor of
$1-c=a+b$.) \ Replacing $a$ by $ps$, $b$ by $qr$ and $c$ by $pr+qs$ correctly
predicts the second part. \ The first part, however, contains an extraneous
term $\ln(ps+qr)$ when expanding the numerator. \ Simulation suggests that our
original formula for $\mathbb{E}\left(  M_{n}\right)  $ is exceedingly
accurate; $\mathbb{E}^{\prime}(M_{n})$ should therefore \emph{not} be employed
in practice.

The associated problem for Geo/Geo/$1$ LAS-DA\ increments remains open, but is
perhaps manageable (owing to their similarity with $Z_{1}$, $Z_{2}$, \ldots,
$Z_{n}$). \ Less feasible, we suspect, would be a rigorous proof of our
asymptotics for two servers or more.

\section{Acknowledgements}

I am thankful to Guy Louchard for introducing me to the Poisson clumping
heuristic (especially recursions for $\nu_{j}$ and $\nu_{-j}$), and to Stephan
Wagner for extracting discrete Gumbel asymptotics in \cite{Fi3-heu} (a
contribution leading to both \cite{Fi0-heu} and the present work). \ Writing
simulation code for M/M/$c$ was facilitated by a theorem in \cite{Sig-heu}
involving order statistics of iid Uniform rvs. \ Umesh Chandra Gupta was so
kind as to send \cite{GG-heu, CG-heu}; Bart Steyaert likewise sent
\cite{GWB-heu}. The creators of R, Julia, Mathematica and Matlab, as well as
administrators of the MIT\ Engaging Cluster, earn my gratitude every day.

\end{document}